\documentclass[a4paper,12pt]{article}

\pdfoutput=1

\usepackage{amsmath}
\usepackage{amsfonts}
\usepackage{amssymb}
\usepackage[latin9]{inputenc}

\numberwithin{equation}{section}

\newtheorem{theorem}{Theorem}[section]
\newtheorem{lemma}{Lemma}[section]
\newtheorem{proposition}{Proposition}[section]
\newtheorem{example}{Example}

\newcommand{\GL}{\Lambda}
\newcommand{\ga}{\alpha}
\newcommand{\gb}{\beta}
\renewcommand{\gg}{\gamma}
\renewcommand{\ge}{\varepsilon}

\newcommand{\GD}{\Delta}
\newcommand{\SA}{{\cal A}}
\newcommand{\tSA}{\widetilde{\cal A}}
\newcommand{\ot}{\otimes}
\newcommand{\gd}{\delta}
\newcommand{\SC}{{\cal C}}

\newcommand{\rd}{{\rm d}}

\newcommand{\SB}{{\cal B}}
\newcommand{\gl}{\lambda}

\newcommand{\SU}{{\cal U}}
\newcommand{\SI}{{\cal I}}
\newcommand{\SW}{{\cal W}}
\newcommand{\kernel}{{\rm kern}\,}

\newcommand{\SV}{{\cal V}}
\newcommand{\SH}{{\cal H}}

\newcommand{\RT}{{\rm T}}
\newcommand{\RS}{{\rm S}}
\newcommand{\rest}{\lceil}
\newcommand{\gs}{\sigma}
\renewcommand{\phi}{\varphi}
\renewcommand{\exp}{{\rm exp}}
\newcommand{\GO}{\Omega}

\renewcommand{\gg}{\gamma}
\newcommand{\beq}{\begin{equation}}
\newcommand{\eeq}{\end{equation}}
\newcommand{\beqo}{\begin{equation*}}
\newcommand{\eeqo}{\end{equation*}}

\newcommand{\id}{{\rm id}}
\newcommand{\bthm}{\begin{theorem}}
\newcommand{\ethm}{\end{theorem}}
\newcommand{\eprop}{\end{proposition}}
\newcommand{\bprop}{\begin{proposition}}
\newcommand{\bexam}{\begin{example}}
\newcommand{\eexam}{\end{example}}
\newcommand{\blem}{\begin{lemma}}
\newcommand{\elem}{\end{lemma}}

\newcommand{\bn}{\par \bigskip \par\noindent}

\newcommand{\pn}{\par\noindent}
\newcommand{\bea}{\begin{eqnarray}}
\newcommand{\eea}{\end{eqnarray}}
\newcommand{\beao}{\begin{eqnarray*}}
\newcommand{\eeao}{\end{eqnarray*}}

\renewcommand{\mit}{ \, \vert \, }
\newcommand{\SN}{{\cal N}}

\newcommand{\SE}{{\cal E}}
\newcommand{\RL}{{\rm L}}
\newcommand{\RM}{{\rm M}}

\newcommand{\ei}{{\bf 1 }}
\newcommand{\op}{\oplus}
\newcommand{\re}{{\rm e}}

\newcommand{\mr}{{\mathbb R}}
\newcommand{\mc}{{\mathbb C}}
\newcommand{\mn}{{\mathbb N}}
\newcommand{\ma}{{\mathbb A}}
\newcommand{\mf}{{\mathbb F}}

\newcommand{\RP}{{\mathbb P}}
\newcommand{\RE}{{\mathbb E}}

\newcommand{\Sqcup}{\amalg}
\newcommand{\sku}{{\cal K}[d]}
\newcommand{\SKU}{{\cal K} \langle d \rangle}
\newcommand{\fn}{\mf _n}
\newcommand{\cfn}{\mc \mf _n}
\newcommand{\exstar}{\exp _{\star}}
\newcommand{\norm}{\vert \vert}
\newcommand{\tR}{\widetilde{R}}
\newcommand{\tC}{\widetilde{C}}
\newcommand{\anz}{\#}
\newcommand{\tphi}{\widetilde{\phi}}
\renewcommand{\Sqcup}{\, \underline{\sqcup} \,}
\renewcommand{\bullet}{\odot}
\newcommand{\sq}{\square}
\newcommand{\tSB}{\widetilde{\SB}}

\newcommand{\tGD}{\widetilde{\GD}}
\newcommand{\bRC}{\mbox{\bf C}}
\newcommand{\bRD}{\mbox{\bf D}}
\newcommand{\bRF}{\mbox{\bf F}}

\begin{document}

\title{Schoenberg Correspondence \\ on Dual Groups \\
{\it \normalsize Dedicated to Wilhelm von Waldenfels}}
\author{ Michael Sch{\"u}rmann \and Stefan Vo{\ss} \\
Department of Mathematics and Computer Science  \\
University of Greifswald  \\
Walther-Rathenau-Stra\ss e 47  \\
17487 Greifswald, Germany   
}
\maketitle
\begin{abstract} 
\noindent
As in the classical case of L\'evy processes on a group,  L\'evy processes on a Voiculescu dual group are constructed from conditionally positive functionals.
It is essential for this construction that Schoenberg correspondence holds for dual groups: The exponential of a conditionally positive functional is a convolution semigroup of states.
\end{abstract}
\section{Introduction}\label{sec:intro}

The original result by Schoenberg \cite{Schoen} says that the Schur exponentials
$(\re  ^{t a_{ij}})_{ij}$, $t \in \mr _+$, of a complex $n \times n$-matrix $A = (a_{ij})_{ij}$
are positive semi-definite iff $A$ is  conditionally positive
semi-definite, i.e. iff
$A ^* = A$ and
$$
\sum_{i, j = 1} ^n \bar z _i \, z_j \, a_{i j} \geq 0
$$
for all complex numbers
$z_1 , \ldots , z_n $ with
$z_1 + \ldots + z_n = 0$.
There are many other examples of Schoenberg correspondence between
conditional positivity and positive semigroups.
An other elementary example is the correspondence between semigroups $P_t =
\re ^{tQ}$
of stochastic matrices and their generator $Q$-matrix which has to have
non-negative
off-diagonal entries and row sums equal to 0.
\par
If $X_t : \GO \to G$, $t \in \mr _+$, is a L\'evy process on a
topological group $G$ 
(that is a $G$-valued stochastic processes with independent
stationary increments) then the distributions $\mu _t$ of $X_t$ form a
convolution semigroup
of probability measures on $G$, i.e.
$
\mu _s \star \mu _t = \mu _{s + t}
$
where the convolution product is defined by
$
(\mu _1 \star \mu _2 ) (f) 
= \int _G \int _G f (xy) \mu _1 (\rd x ) \mu _2 (\rd y ) 
$
for probability measures $\mu$, $\nu$ and bounded continuous functions $f$ on
$G$.
If $\mu _t$ is weakly continuous, one can define the generator of the
convolution
semigroup on an appropriate $*$-algebra of functions on $G$.
This leads to the L\'evy-Khintchine formula in the case of $G = \mr ^d$, or,
more generally,
to Hunt's formula if $G$ is a Lie group.
Again the generator $\psi$ is conditionally positve,
in the sense that $\psi$ is hermitian and $\psi (f) \geq 0$ for functions
$f \geq 0$ vanishing at the unit element of $G$; see \cite{Hunt,vW2}.
In the case when $G$ is a locally compact abelian group or a compact group, one
can choose as space of functions on $G$ the $*$-algebra formed by the
coefficient 
functions of continuous irreducible representations of $G$.
This coefficient algebra $\SB$  is a Hopf $*$-algebra, and convolution semigroups of
states on $G$ are precisely
given by the conditionally positive linear functionals on the coefficient algebra
where conditionally
positive now means hermitian and
$$
\psi (f) \geq 0 \ \mbox{for} \ f \geq 0 , f \in \kernel \gd .
$$
The functional $\gd$ is the counit of the coefficient algebra $\SB$ and convolution of linear functionals 
$\phi _1$ and $\phi _2$ on  $\SB$ is given by
\begin{eqnarray}
\label{conv1}
\phi _1 \star \phi _2 = (\phi _1 \ot \phi _2 ) \circ \GD ,
\end{eqnarray}
where $\GD$ is the comultiplication of $\SB$.
\par
We describe the  mechanism  of constructing L\'evy processes.
Starting from a conditionally positive linear functional $\psi$ on $\SB$,
we obtain a convolution semigroup $\phi _t$ of linear functionals on $\SB$ as the convolution exponentials
$\re _{\star} ^{t \psi}$ of $\psi$.
Now it is important that Schoenberg correspondence holds which means
that $\phi _t$ are positive so that the convolution semigroup consists of states which again are in
1-1-correspondence to probability measures on $G$.
The convolution semigroup defines a projective system of finite-dimensional distributions which
by Kolmogorov's theorem allows to construct a L\'evy process on $G$ whose convolution semigroup
is given by $\phi _t$.
This establishes, up to stochastic equivalence of stochastic processes,
a 1-1-correspondence between conditionally positive linear functionals
on the coefficient algebra of $G$ and L\'evy processes on $G$.
\par
The Hopf $*$-algebras arising from locally compact abelian or compact groups
are algebras of functions and as such are commutative.
If one generalizes to arbitrary Hopf $*$-algebras (for instance, compact quantum groups),
a notion of noncommutative (quantum) L\'evy processes has been introduced; see \cite{MSchue93}.
These quantum L\'evy processes are again given by conditionally
positive linear functionals, now on the Hopf $*$-algebra $\SB$ where conditionally positive now means hermitian and
\begin{eqnarray}
\psi (b ^* b ) \geq 0 \ \mbox{for} \ b \in \kernel \gd .
\end{eqnarray}
The increments of quantum L\'evy processes on Hopf $*$-algebras are independent in the sense of {\it noncommutative tensor independence} where
two sub-algebras are called tensor independent if they commute and if expectations factorize.
\par
The tensor product of linear functionals is closely related to the classical notion of stochastic independence.
Two random variables $X_1 : \GO  \to G_1$, $X_2 : \GO  \to G_2$ are independent if their joint distribution is the tensor product of the marginal distributions that is if 
\begin{eqnarray}
\label{indep1}
\RP _{(X_1 , X_2 )} = \RP _ {X_1} \ot \RP _{X_2} .
\end{eqnarray}
Identify the underlying probability measure $\RP$ with its expectation
$\RE (F) = \int  F \rd P$, $F \in \RL ^{\infty}(\GO )$,
which is a positive normalized linear functional on the $*$-algebra $\RL ^{\infty} (\GO )$.
If we think of $\RP _{X_1}$ and $\RP _{X_2}$ as their expectations $\RE _{X_1}$ and $\RE _{X_2}$ on the $*$-algebras
$\RL ^{\infty} (G_1 )$ and
$\RL ^{\infty} (G_2 )$ of
bounded measurable functions on $G_1$ and on $G_2$ respectively, and of the joint distribution as the expectation
$\RE _{(X_1, X_2 )}$ on the tensor product
$
\RL ^{\infty} (G_1 \times G_2 ) = \RL ^{\infty} (G_1 ) \ot \RL ^{\infty} (G_2 )
$ 
then (\ref{indep1}) becomes
\begin{eqnarray}
\label{indep2}
\RE _{(X_1, X_2 )} = \RE _{X_1} \ot \RE _{X_2}.
\end{eqnarray}
Define the $*$-algebra homomorphisms $j_1 : \RL ^{\infty} (G_1 ) \to \RL ^{\infty}(\GO ) $ and
$ \RL ^ {\infty}(G_2 )  \to \RL ^{\infty}(\GO ) $ by $j_1 (f_1 ) = f_1 \circ X_1$
and $j_2 (f_2 ) = f_2 \circ X_2$ and introduce the $*$-algebra homomorphism
$$
j_1 \ot j_2 : \RL ^{\infty} (G_1 ) \ot \RL ^{\infty} (G_2 ) \to \RL ^{\infty} (\GO ) 
$$
by
$$
(j_1 \ot j_2 ) (f_1 \ot f_2 ) = j_1 (f_1 )\, j_2 (f_2 ) .
$$
Then $\RE \circ (j_1 \ot j_2 )$ is the joint distribution of $X_1$
and $X_2$, and (\ref{indep2}) reads
\begin{eqnarray}
\label{indep3}
\RE \circ (j_1 \ot j_2 ) = (\RE \circ j_1 ) \ot (\RE \circ j_2 )
\end{eqnarray}
\par
It is remarkable that in a noncommutative world there is more than one
possibility for a notion of independence.
One of these notions, called tensor independence, as was pointed out above is closely related to classical independence and is the noncommutative independence chosen for L\'evy processes on Hopf $*$-algebras.
In his papers on the broadening of spectral lines \cite{vW1} W. von Waldenfels used another notion of independence, in
some respect the most simple, which now is called {\it Boolean} independence, because Boolean lattices appear when moments are calculated from their cumulants. 
\par
A central role is played by {\it free} independence or {\it freeness} which was introduced
by D. Voiculescu \cite{Voi1}.
N. Muraki showed that under certain natural axioms there are exactly five notions of noncommutative independence (\cite{MurN02,MurN03}, see also \cite{Spe,BGhSc02}), namely tensor, Boolean, free and monotonic and anti-monotonic independence.
\par
The joint distribution of two classical random variables lives on the tensor product which is the commutative algebra
\lq freely\rq \ generated by the algebras $\RL ^{\infty} (G_1 )$ and $\RL ^{\infty} (G_2 )$.
In the noncommutative case the tensor product of algebras has to be replaced by the free product of algebras.
The classification result of Muraki classifies the \lq natural\rq \ products of linear functionals which assign to each pair of algebras $(\SB _1 , \SB _2 )$ and linear functionals $(\phi _1 , \phi _2 )$, $\phi _1 : \SB _1 \to \mc$,
$\phi _2 : \SB _2 \to \mc$ a linear functional $\phi _1 \bullet \phi _2 : \SB _1 \sqcup \SB _2 \to \mc$.
The product  $\bullet$ replaces the tensor product of (\ref{indep2}) of the classical case.
If we understand by a noncommutative (quantum) probability space a $*$-algebra $\SA$ equipped with a state
$\RE : \SA \to \mc$ and by a quantum random variable on a $*$-algebra $\SB$ over a
quantum probability space $(\SA , \RE )$ a homomorphism $j: \SB \to \SA$ of $*$-algebras, two random variables
$j_1 : \SB _1 \to \SA$ and $j_2 : \SB _2 \to \SA$ are called independent
(with respect to an independence given by a Muraki natural product $\bullet$) if
\begin{eqnarray}
\label{indep4}
\RE \circ (j_1 \sqcup j_2 ) = (\RE \circ j_1 ) \bullet (\RE \circ j_ 2)
\end{eqnarray}
which precisely is the noncommutative version of (\ref{indep3}).
\par
In this paper, general quantum L\'evy processes are considered, the independence of increments coming from one of the five notions of independence  of Muraki's classification.
To treat the general independence case, Hopf $*$-algebras have to be replaced by their
\lq free\rq \ counterparts where tensor products of algebras are replaced by free products of algebras.
Such objects appeared already in work of D. Voiculescu \cite{Voi87} and had been called \lq dual groups\rq .
We will call \lq dual semigroup\rq \ a $*$-algebra $\SB$ equipped with a comultiplication $\GD$ which is
a $*$-algebra homomorphism from $\SB$ to the free  product $\SB \sqcup \SB$ of $\SB$ with itself such that 
coassociativity and the counit property hold.
If there is also an antipode we will speak of (algebraic) dual groups.
\par
Suppose that $(\SB, \GD )$ is a dual semigroup.
In addition let there be given a (fixed) natural product $\bullet$ with its associated
notion of noncommutative independence; see \cite{MurN03,BGhSc02}.
We define the convolution product of two linear functionals
$\phi _1$, $\phi _2$ on $\SB$ by
\begin{equation}
\label{conv2}
\phi _1 \star \phi _2 = (\phi _1 \bullet \phi _2 ) \circ \GD 
\end{equation}
in complete analogy to the tensor case (\ref{conv1}).
Quantum L\'evy processes on the dual semigroup $\SB$ are again determined by convolution semigroups
of states on $\SB$.
Convolution exponentials $\re _{\star} ^{ \psi }$ of linear functionals $\psi$ on $\SB$ can be defined as before (see Section \ref{sec:schoen}).
Schoenberg correspondence (Theorem \ref{theorem}) says that the (point-wise) continuous convolution semigroups
$\phi _t$ of states on $\SB$ are precisely given by $\phi _t = \re _{\star} ^{t \psi}$ with
$\psi$ the (point-wise) derivative of $\phi _t$ and $\psi$  conditionally positive.
\par
In Section \ref{sec:levy} we start from a conditionally positive linear functional $\psi$ on
a dual semigroup $\SB$.
Using Schoenberg correspondence, we associate with it a convolution semigroup of states on $\SB$.
By an inductive limit procedure we then construct a quantum L\'evy process on $\SB$ with 
convolution semigroup given by $\phi _t = \re _{\star} ^{t \psi}$.

\section{Preliminaries}\label{sec:pre}

Algebras will be over the complex numbers and will assumed to be associative.
A $*$-algebra is an algebra $\SA$  with an involution $*$, i.e. an anti-linear map
$a \mapsto a ^{*}$ on $\SA$ such that $(ab) ^{*} = b ^* a ^*$ and $(a ^* )^* = a$.
A unital algebra is an algebra such that there exists an element $\ei$ (called the unit element) in $\SA$ with
$a \, \ei = a = \ei \, a$.
A (counital) coalgebra is a triplet $(\SC , \GD , \gd )$ consisting of a (complex) vector space
$\SC$ and linear mappings $\GD : \SC \to \SC \ot \SC$ and $\gd : \SC \to \mc$ such that
$(\GD \ot \id) \circ \GD =  (\id \ot \GD ) \circ \GD$ and
$(\gd \ot \id ) \circ \GD = \id = (\id \ot \gd ) \circ \GD$
where for vector spaces $\SV$ and $\SW$ we write $\SV \ot \SW$ for the vector space tensor product.
Either $\SC$ is the trivial vector space or there exists an element $e \in \SC$ with
$\gd e = 1$. In the latter case $\SC = \mc e \op \SC_0$ with $\SC _0 = \kernel \gd$, and for
$c \in \SC$, $c _0 = c - \gd (c) e$, we have
$\GD c - (\gd (c) e \ot e + e \ot c _0 + c _0 \ot e ) \in \SC _0 \ot \SC _0$.
In particular, $\GD e = e \ot e + B$, $B \in \SC _0 \ot \SC _0$, and $\GD c = e \ot c + c \ot e + B$, $B \in \SC _0 \ot \SC _0$, for $c \in \SC _0$.
\par
A bialgebra is a coalgebra $(\tSB , \tGD , \gd )$ where $\tSB$
is a unital algebra such that $\tGD$, $\gd$ are algebra homomorphisms.
If we put $\SB = \kernel \gd$, then
$\tGD \SB \subset \SB \oplus \SB \oplus (\SB \ot \SB ) = : \SB \ot _0 \SB$, and the pair $(\SB , \GD )$, $\GD = \tGD \rest \SB : \SB \to \SB \ot _0 \SB$ consists of an algebra $\SB$ and a \lq comultiplication\rq \ $\GD$ such that
\begin{equation}
\label{coass0}
(\GD \ot _0 \id ) \circ \GD = (\id \ot _0 \GD ) \circ \GD
\end{equation}
and
\begin{equation}
\label{counit0}
(0 \ot _0 \id ) \circ \GD = \id = (\id \ot _0 0 ) \circ \GD .
\end{equation}
This shows that a bialgebra equivalently can be defined to be a pair $(\SB , \GD )$ consisting of an algebra $\SB$ and an algebra homomorphism $\GD : \SB \to \SB \ot _0 \SB$ such that \eqref{coass0} and \eqref{counit0} hold.
\par
For an index set $I \neq \emptyset$, we put
$$
\ma (I) =
\{ (\ge _1 , \ldots , \ge _n) \mit n \in \mn, \ge _l \in I, l = 1, \ldots , n, \ge _l \neq \ge _{l + 1},
l = 1, \ldots , n - 1 \} .
$$
For a family $(\SA _i ) _{i \in I}$ of algebras and $\ge = (\ge _1 , \ldots , \ge _n ) \in \ma (I)$ we denote by
$\SA _{\ge}$ the algebraic tensor product
$\SA _{\ge} = \SA _{\ge _1} \ot \ldots \ot \SA _{\ge _n}$ of the algebras
$\SA _{\ge _1}, \ldots , \SA _{\ge _n}$.
Define the free product $\bigsqcup _{i \in I} \SA _i$ of the family $(\SA _i ) _{i \in I}$
as the vector space direct sum
$$
\bigsqcup _{i \in I} \SA _i = \bigoplus _{\ge \in \ma (I)} \SA _{\ge}
$$
with multiplication given by
\beao
&& (a_1 \ot \ldots \ot a_n ) \, (b_1 \ot \ldots \ot b_m )  \\
&& \qquad = \left\{ \begin{array}{lll}
a_1 \ot \ldots \ot a_n \ot b_1 \ot \ldots \ot b_m & \mbox{if} & \ge _n \neq \gamma _m  \\
a_1 \ot \ldots \ot a_{n - 1} \ot (a_n b_1 ) \ot b_2 \ot \ldots \ot b_m & \mbox{if} & \ge _n = \gamma _m 
\end{array}
\right.
\eeao
for $n, m \in \mn$, $\ge = (\ge _1 , \ldots , \ge _n ) \in \ma (I)$,
$\gg = (\gamma _1 , \ldots , \gamma _m ) \in \ma (I)$ and
$a_1 \ot \ldots \ot a_n \in \SA _{\ge}$, $b_1 \ot \ldots \ot b_m \in \SA _{\gamma}$.
For example, if $I = \{ 1, 2 \}$
$$
\SA _1 \sqcup \SA _2 = \SA _1 \op \SA _2 \op (\SA _1 \ot \SA _2) \op (\SA _2 \ot \SA _1) \op \ldots
$$
The free product is the co-product in the category of algebras, i.e.
given two algebra homomorphisms $j_1 : \SB _1 \to \SA $ and $j_2 : \SB _2 \to \SA$ with the same target $\SA$, there is a unique algebra homomorphism $j_1 \sqcup j_2 : \SB _1 \sqcup \SB _2 \to \SA $ such that
$j_{1/2} = (j_1 \sqcup j_2 ) \circ i_{1/2}$
where $i_1$, $i_2$ denote the natural embeddings of $\SB _1$, $\SB _2$ into $\SB _1 \sqcup \SB _2$.
We frequently write $j_1 \Sqcup j_2$ for $(i_1 \circ j_1 ) \sqcup (i_2 \circ j_2 ) : \SB _1 \sqcup \SB _2 \to \SA _1 \sqcup \SA _2$ for algebra homomorphisms $j_{1/2} : \SB _{1/2} \to \SA _{1/2}$.
The free product $\SA _1 \sqcup _{\ei} \SA _2$ of unital algebras $\SA_1 , \SA _2$ is obtained from $\SA _1 \sqcup \SA _2$ by dividing by the ideal generated by $\ei _{\SA _1} - \ei _{\SA _2}$.
Then $\sqcup _{\ei}$ is the co-product in the category of unital algebras.
\par
We follow \cite{Fra06} and define (stochastic) independence in the language of category theory.
Let $(\bRC , \square , i )$ be a {\it tensor category with injections}, i.e. $(\bRC , \square )$ is a tensor category such that for each pair  $ C_1, C_2 $ of objects there exists a pair $i_{C_1} , i_{C_2}$ of morphisms $i_{C_1} : C_1 \to C_1 \square C_2$, $i_{C_2} : C_2 \to C_1 \sq C_2$, such that for any pair $j_1 : C_1 \to D_1 $, $j_2 : C_2 \to D_2$ of morphisms we have
\begin{align*}
(j_1 \sq j_2 ) \circ i _{C_1} &= i_{D_1} \circ j_1  , \\
(j_1 \sq j_2 ) \circ i _{C_2} &= i_{D_2} \circ j_2 .
\end{align*}
Then two morphisms $j_1 : C_1 \to C$, $j_2 : C_2 \to C$ with the same target are called {\it independent} if there exists a morphism $j: C_1 \sq C_2 \to C$ such that $j_1 = j \circ i _{C_1}$ and $j_2 = j \circ i_{C_2}$.
For example, consider the category formed by \lq dual probability spaces\rq \ that is by pairs $(C , \phi )$ with $C$ a commutative von-Neumann algebra and $\phi$ a normal state on $C$. 
This is a tensor category with injections if we choose the von-Neumann algebra tensor product with the tensor product of states and the natural injections. The morphisms from $(C, \phi )$ to $(D , \theta )$ are the von-Neumann algebra homomorphisms $j$ with $\phi = j \circ \theta$, that is they are precisely the random variables.
Two random variables are stochastically independent in the classical sense iff they are independent in the above sense of categories with injections.
\par
Since we will stay in an algebraic framework, we consider tensor products in the category formed by pairs $(\SB , \phi )$ where $\SB$ is an algebra and $\phi : \SB \to \mc$ is a linear functional on $\SB$.
 There is a type of tensor product with injections in this category given by 
\[
(\SB _1 , \phi _1 ) \sq (\SB _2 , \phi _2 ) =
(\SB _1 \sqcup \SB _2 , \phi _1 \odot \phi _2 )
\]
where $\phi _1 \odot \phi _2$ is a linear functional on the free product $\SB _1 \sqcup \SB _2$ of algebras such that the product $\odot$ satisfies the axioms
\begin{gather*}
 (\phi _1 \odot \phi _2 ) \circ i _{1/2 } = \phi _{1/2}     \tag{A1}\label{A1}
\\ 
(\phi _1 \odot \phi _2 ) \odot \phi _3 = \phi _1 \odot (\phi _2 \odot \phi _3 )    \tag{A2}\label{A2}
\\
(\phi _1 \circ j_1 ) \odot (\phi _2 \circ j_2 ) =
 (\phi _1  \odot \phi _2 )\circ (j_1 \Sqcup j_2 )  \tag{A3}\label{A3}
\end{gather*}
Consider the additional axioms
\beq
(\phi _1 \odot \phi _2 ) (b_1 b_2) = (\phi _1 \odot \phi _2 ) (b_2 b_1) =
\phi _1 (b_1 ) \phi _2 (b_2 )
\tag{A4}\label{A4}
\eeq
for all $b_1 \in \SB _1 , b_2 \in \SB _2$, and
\beq
\phi _1 \odot \phi _2 = \phi _2 \odot \phi _1 . \tag{A5}\label{A5}
\eeq
N. Muraki \cite{MurN02,MurN03} showed that there are exactly five products 
satisfying  \eqref{A1}-\eqref{A4}, the tensor product, the free product \cite{Voi1}, the Boolean product \cite{vW1}, and the monotonic and anti-monotonic products \cite{MurN01,Lu}.
It was shown in \cite{BGhSc02,Spe} that the tensor, the free and the Boolean products are the only three products satisfying \eqref{A1}-\eqref{A5}.
An independence coming from a product with \eqref{A1}-\eqref{A3} will be called a {\it $\odot$-independence}.
\par
 A linear functional $\phi$ on a $*$-algebra $\SA$ is called  hermitian if  $\phi (a ^* ) = \overline{\phi (a)}$  for all $a \in \SA$. We call $\phi$ {\it conditionally positive} if $\phi$ is hermitian and if $\phi (a ^* a ) \geq 0$ for all $a \in \SA$. In this paper, we call $\phi$ a {\it state} if we have $\tphi (a^* a ) \geq 0$ for all $a \in \tSA$ where $\tphi : \tSA \to \mc$ is the normalized linear extension of $\phi$ to $\tSA = \mc \ei \oplus \SA$.
 A state is conditionally positive whereas the converse is not always true.
For example, on the $*$-algebra of complex polynomials in one self-adjoint indeterminate $x$, with constant part equal to 0 the linear functional $\psi$ with $\psi (x ^n ) = \gd _{2, n}$ is conditionally positive but not a state.
The point-wise limit of states is a state. If $\SI$ is a two-sided $*$-ideal of the $*$-algebra $\SA$, then for a state $\phi$ on $\SA$ which vanishes on $\SI$ we have that $\hat \phi$, $\hat \phi (a + \SI  ) = \phi (a)$, is a state on $\SA /\SI$.
\par
We say that a $\odot$-independence is {\it positive} if for two $*$-algebras $\SA _1$ and $\SA _2$ and states $\phi _1$ and $\phi _2$ on $\SA _1$ and $\SA _2$ respectively, the product $\phi _1 \bullet \phi _2$ is a state.
Another notion of states is the following. Call $\phi$ a strong state if $\min \{ \gl \in \mc \mit \phi _{\gl} (a ^* a ) \geq 0 \ \forall a \in \tSA \} = 1$
where $\phi _{\gl}$ is the extension of $\phi$ to $\mc \ei \oplus \SA$ with $\phi _{\gl} (\ei ) = \gl$. 
Then a strong state is a state, and the converse is false in general. We say that a  $\odot$-independence is strongly positive if the product of two strong states is a strong state. 
Then each positive $\odot$-independence is strongly positive.
This holds because $\phi _1$ and $\phi _2$ are the restrictions of $\phi _1 \bullet \phi _2$ to $\SA _1$ and $\SA _2$.
It is well-known that Muraki's five notions of independence are positive.
In fact,
\bprop
The only positive $\odot$-independences are Muraki's five. In particular,
a $\odot$-independence is strongly positive iff it is positive.
\eprop
\pn
{\it Proof}:
We show that a strongly positive $\odot$-independence must be one of Muraki's five. Take $\SB _1 = \SB _2 = \mc [x]$ and $\phi (x^n ) = 1$, $n \in \mn$. 
Then $\phi$ is a strong state. If $\odot$ is strongly positive $\phi \odot \phi$ must be a state.
It follows from \cite{BGhSc02}, Lemma 2.1, that with
$\mc [x ] \sqcup \mc [x] = \mc \langle x, y \rangle$ we must have $(\phi  \odot \phi )(xy) = q_1 $ and $\phi  \odot \phi (yx) = q_2 $ for some complex constants $q_1$ and $q_2$. 
We have
$$
0 \leq (\phi \odot \phi )\bigl( (\gl \ei + \ga x + \gb y )^* (\gl \ei + \ga x + \gb y )\bigr)
$$ 
for all $\gl, \ga , \gb \in \mc$ and the matrix
$$
\left( \begin{array}{ccc} 1 & 1 & 1  \\
1 & 1 & q_1  \\
1 & q_2 & 1 \end{array} \right)
$$ 
must be positive semi-definite which forces
 $q_1 = q_2 = 1$. It follows that a strongly positive $\odot$-independence satisfies (A4) so that Muraki's result can be applied.$\square$
\par
A {\it dual semigroup} is a pair $(\SB , \GD )$ consisting of a $*$-algebra $\SB$ and a $*$-algebra homomorphism $\GD : \SB \to \SB \sqcup \SB$ such that
\[
(\GD \Sqcup \id ) \circ \GD = (\id \Sqcup \GD ) \circ \GD
\]
and
\[
(0 \sqcup \id ) \circ \GD = \id = (\id \sqcup 0 ) \circ \GD ;
\]
cf. \cite{Zhang,Voi87} and  \cite{BGhSc}.
Put $\tSB = \mc \ei \oplus \SB$, $\tGD : \tSB \to \tSB \sqcup _{\ei} \tSB
= \mc \ei \oplus \SB \sqcup \SB$,
$\tGD \rest \SB = \GD$, $\tGD \ei = \ei$, $\gd : \tSB \to \mc$, $\gd \rest \SB = 0$, $\gd \ei = 1$. 
Then the triplet $(\tSB , \tGD , \gd )$ satisfies
\begin{align}
\label{coass}
(\tGD \Sqcup _{\ei} \id ) \circ \tGD &= (\id \Sqcup _{\ei} \tGD ) \circ \tGD
\\ 
\label{counit}
(\gd \Sqcup _{\ei} \id ) \circ \tGD &= \id = (\id \Sqcup _{\ei} \gd ) \circ \tGD .
\end{align}
Conversely, given a triplet $(\tSB , \tGD , \gd )$ such that $\tSB$ is a unital $*$-algebra and $\tGD : \tSB \to \tSB \sqcup _{\ei} \tSB$, $\gd : \tSB \to \mc$ are unital $*$-algebra homomorphisms with \eqref{coass} and \eqref{counit}, it can be shown  that the pair $(\kernel \gd , \tGD \rest \kernel \gd )$ is a dual semigroup; see \cite{BGhSc}.
A dual semigroup is called a {\it dual group} if there is an endomorphism $S$ on $\SB$ such that $(S \sqcup \id ) \circ \GD = 0 = (\id \sqcup S ) \circ \GD$.
\par
Let $\bRC$ and $\bRD$ be two categories and let $\bRF$ be a functor from $\bRC$ to $\bRD$. For an object $D$ in $\bRD$ a {\it universal pair} or arrow
\cite{MacLane} from $D$ to $\bRF$ is a pair $(C , i )$ with $C$ an object in $\bRC$ and $i$ a morphism $i : D \to \bRF (C )$ such that the following universal property is fulfilled. For each object $A$ in $\bRC$ and morphism $k : D \to \bRF (A)$ there is a unique morphism $j : C \to A$ such that $k = \bRF (j ) \circ i$.
In the case when $\bRC$ is the category of algebras, $\bRD$ is the category of vector spaces, and $\bRF$ is the forgetful functor, a universal pair from a vector space $\SV$ to $\bRF$ can be realized as the tensor algebra over $\SV$ which is the vector space direct sum of the vector space tensor powers $\SV ^{\ot n}$ of $\SV$. This is an algebra with multiplication $(v_1 \ot \ldots v_n ) (w_1 \ot \ldots \ot w_m ) = v_1 \ot \ldots \ot v_n \ot w_1 \ot \ldots \ot w_m$, $v_1, \ldots, v_n , w_1 , \ldots , w_m \in \SV$. The morphism $i_{\SV}$ is given by the natural embedding.
The unique morphism $j$ associated with a morphism $k$ is denoted by $\RT (k)$.
Similarly, the tensor $*$-algebra $\RT (\SV )$ over a $*$-vector space $\SV$, that is  a vector space  $\SV$ equipped with an anti-linear self-inverse map $v \mapsto v^*$, is the universal pair from $\SV$ to the forgetful functor from the category of $*$-algebras to the category of $*$-vector spaces. The involution of $\RT (\SV )$ is given by $(v_1 \ot \ldots \ot v_n  ) ^* = v_n ^* \ot \ldots \ot v_1 ^*$.
\par
Now let $\bRC$ be the category of {\it commutative} algebras and let $\bRD$ again be  the category of vector spaces with the forgetful functor from $\bRC$ to $\bRD$. The universal pair is denoted by $(\RS (\SV ) , i_{\SV})$, and a realization of $\RS (\SV )$ is the symmetric tensor algebra over $\SV$ which is the quotient of $\RT (\SV )$ by the ideal generated by $v \ot w - w \ot v$, $v, w \in \SV$. We frequently write $v_1 \ot _s \ldots \ot _s v_n$ for the equivalence class of $v_1 \ot \ldots \ot v_n$.
Of course, there is also the symmetric tensor $*$-algebra over a $*$-vector space $\SV$.
If $\{ v_i \mit i \in I \}$ is a (self-adjoint) vector space basis of $\SV$, then $\RS (\SV )$ can be identified with the polynomial algebra $\mc [ x_i ; i \in I ]$.
\par
It follows from \cite{BGhSc02}, Lemma 2.1, that  if $\odot$ satisfies \eqref{A3} there are linear maps
\[
\gs _{\SB _1 , \SB _2} : \SB _1 \sqcup \SB _2 \to 
\RS (\SB _1 ) \ot \RS (\SB _2 ) \cong \RS (\SB _1 \oplus \SB _2 )
\]
such that
\[
\phi _1 \odot \phi _2 = ( \RS (\phi _1 ) \ot \RS (\phi _2 )) \circ \gs _{\SB _1 , \SB _2} .
\]
By Theorem 3.4 of \cite{BGhSc}, for a fixed $\odot$-independence and a dual semigroup $(\SB , \GD )$, we can form the commutative $*$-bialgebra $(\RS (\SB ),  \RS (\gs \circ \GD ))$ where we put $\gs = \gs _{\SB , \SB}$. 
Thus a $\odot$-independence gives rise to a functor from the category of dual semigroups to the category of commutative $*$-bialgebras.
Now put, for a fixed $\odot$-independence,
\begin{equation}
\label{conv}
\phi _1 \star \phi _2 = (\phi _1 \odot \phi _2 ) \circ \GD
\end{equation}
for linear functionals $\phi _1$ and $\phi _2$ on a dual semigroup $\SB$.
Then (see \cite{BGhSc02}) 
\begin{equation}
\label{hom}
\RS (\phi _1 \star \phi _2 ) =
 \RS (\phi _1 ) \star \RS (\phi _2 ) 
\end{equation}
 where the second convolution product is with respect to the comultiplication $\RS (\gs \circ \GD )$.

\section{Schoenberg correspondence}\label{sec:schoen}

The intersection of two coalgebras is again a coalgebra so that the sub-coalgebra generated by a subset of a coalgebra is well-defined. The fundamental theorem of coalgebras (see e.g. \cite{DNR}) says that the sub-coalgebra generated by a single element (and thus by a finite number of elements) is finite-dimensional. It follows that a coalgebra is the inductive limit of its finite-dimensional sub-coalgebras. For a linear functional $\psi$ on a coalgebra $\SC$ define the linear map $T_{\psi} $ on $\SC$ by $T_{\psi} = (\id \ot \psi ) \circ \GD$.
Then $T$ is a unital algebra homomorphism from the convolution algebra formed by linear functionals on $\SC$ to the algebra of linear operators on $\SC$ and $L \mapsto \gd \circ L$ is the left inverse of $T$.
Moreover, $T _{\psi}$ leaves invariant all sub-coalgebras of $\SC$.
Denote by $\re ^{T _{\psi}}$ the inductive limit of the (matrix) exponentials of the restrictions of $T_{\psi}$ to finite-dimensional sub-colagebras.
Put $\exp _{\star} \psi : = \gd \circ \re ^{T_{\psi}}$.
It follows that the series 
\begin{equation}
\label{convexp}
\sum_{n = 0}^\infty  \frac{\psi ^{\star n}}{n !} (c)
\end{equation}
converges for all $c \in \SC$ and that this limit equals $\exp _{\star} \psi$. We have
$$
\exstar (\psi _1 + \psi _2 ) = (\exstar \psi _1 ) \star (\exstar \psi _2 )
$$
if $\psi _1 \star \psi _2 = \psi _2 \star \psi _1$ and
$$
(\exstar \psi )(c) = \lim_{n \to \infty}
(\gd + \frac{\psi}{n} ) ^{\star n} (c) .
$$
More generally, \cite{ScVo,ScSkVol}
\blem
\label{lemma1}
Let $\psi$ be a linear functional on a coalgebra $\SC$.
Suppose that $R_n$, $n \in \mn$, are linear functionals on $\SC$ such that for each $b \in \SC$ there is a constant $C_b \in \mr_+$ with
\begin{equation}
\label{inequ}
\vert R_n (b) \vert \leq \frac{1}{n ^2} \, C_b  \  \  \forall n \in \mn .
\end{equation}
Then
$$
(\gd + \frac{\psi}{n} + R_n ) ^{\star n}
$$
converges to $\exstar \psi$ point-wise.
\elem
\pn
{\it Proof}:
By the fundamental theorem of coalgebras we can assume that $\SC$ is finite-dimensional.
Then with some norm $\norm \ \norm$ on $\SC$
\beao
\norm T_{R_n} \norm &=& \norm (\id \ot R_n ) \circ \GD \norm  \\
&\leq&
\norm \id \ot R_n \norm \, \norm \GD \norm   \\
&=& \norm R_n \norm \, \norm \GD \norm .
\eeao
Choose a vector space basis $\{ b_1 , \ldots , b_k \}$ of
$\SC$.
Then for $\ga _1 , \ldots , \ga _n \in \mc$
\beao
&&\vert R_n (\ga _1 b_1 + \ldots + \ga _k b_k ) \vert  \\
&&\qquad \leq \max_{1 \leq j \leq k} \vert R_n (b_j ) \vert \,
(\vert \ga _1 \vert + \ldots + \vert \ga _k \vert )   \\
&&\qquad  \leq  \frac{1}{n^2 } (\max_{1 \leq j \leq k} C_{b_j} )
(\vert \ga _1 \vert + \ldots + \vert \ga _k \vert )
\eeao
which implies $\norm T_{R_n} \norm \leq \frac{1}{n^2} \, C$ for some constant $C$. Now
$$
(\id + \frac{T_{\psi}}{n} + T_{R_n} ) ^n \to \re ^{T_{\psi}}
$$
and thus 
$$
(\gd + \frac{\psi}{n} + R_n ) ^{\star n}  \to  \re ^{\psi} . \square
$$
A family $(\phi _t )_{t\in \mr _+}$ of linear functionals on a coalgebra $\SC$ is called a {\it continuous convolution semigroup} (CCSG) if
$\phi _{s + t} = \phi _s \star \phi _t$, $\phi _0 = \gd$, and 
$\phi _t \to \gd$ point-wise for $t \to 0+$. 
For a CCSG the operators $T_{\phi _t}$ form a semigroup of linear operators on $\SC$.
Using the fundamental theorem of coalgebras and a well-known result for continuous semigroups of matrices, we obtain that 
$\lim_{t \to 0+} \frac1t (\phi _t - \gd )$ exists point-wise and that we have
$\exstar (t \psi ) = \phi _t$  for the limiting functional $\psi$.
It follows that the CCSGs are exactly given by the convolution exponentials $\exstar (t \psi )$.
\par
Now let $\psi$ be a linear functional on a dual semigroup $\SB$.
Moreover, fix a $\odot$-independence, and put $D (\psi ) (B) =
\frac{\rd}{\rd \, t} \RS (t \, \psi )(B)  \vert_{t = 0}$ for $B \in \RS (\SB )$.
We define $\exstar \psi $ point-wise by $(\exstar D (\psi ) ) \circ i _{\SB}$.
Then $\exstar (t \psi)$, $t \geq 0$, form a CCSG of linear functionals on $\SB$ with the convolution product now given by \eqref{conv}; see \cite{BGhSc02}.
We have \cite{ScVo}
\blem
\label{lemma2}
Let $\psi$ be a linear functional on a dual semigroup $\SB$.
Suppose that $R_n$, $n \in \mn$, are linear functionals on $\SB$ such that for each $b \in \SB$ there is a constant $C_b \in \mr_+$ with
$$
\vert R_n (b) \vert \leq \frac{1}{n ^2} \, C_b  \  \  \forall n \in \mn .
$$
Then
$$
(\frac{\psi}{n} + R_n ) ^{\star n}
$$
converges to $\exstar \psi$ point-wise.
\elem
\pn
{\it Proof}:
By \eqref{hom} we have
$$
S(( \frac{\psi}{n} + R_n ) ^{\star n} ) = S (\frac{\psi}{n} + R_n ) ^{\star n}
$$
and 
\beao
(\frac{\psi}{n} + R_n ) ^{\star n} &=&
S((\frac{\psi}{n} + R_n ) ^{\star n} ) \rest \SB   \\
&=&
S(\frac{\psi}{n} + R_n ) ^{\star n } \rest \SB .
\eeao
Moreover, $\exstar \psi = \exstar D(\psi ) \rest \SB$.
We will prove that 
\[
S( \frac{\psi}{n} + R_n ) ^{\star n} \to \exstar D(\psi ).
\]
For $b_1 , \ldots , b_k \in \SB$, $k \geq 1 $,
\beao
&& S(\frac{\psi}{n} + R_n ) (b_1 \ot _s \ldots \ot _s b_k )   \\
&& \qquad = (\frac{\psi}{n} + R_n ) (b_1 ) \ \ldots (\frac{\psi}{n} + R_n ) (b_k )  \\
&& \qquad = \sum_{A \subset \{ 1 , \ldots , k \}} \frac{1}{n^{\anz A}} \prod_{j \in A} \psi (b_j ) \prod_{j \notin A} R_n (b_j )   \\
&& \qquad = R_n (b_1 ) \ldots R_n (b_k ) + \frac1n \bigl( \psi (b_1 ) R_n (b_2 ) \ldots R_n (b_k )   \\
&&\qquad \quad + R_n (b_1 ) \psi (b_2 ) R_n (b_3 ) \ldots R_n (b_k ) + \ldots + R_n (b_1 ) \ldots R_n (b_{k - 1} ) \psi (b_k ) \bigr)  \\
&& \qquad\qquad + \frac{1}{n^2} T_n (b_1 \ot _s \ldots \ot _s b_k )
\eeao
with $\vert T_n (b_1 \ot _s \ldots \ot _s b_k ) \vert \leq D_1 $ for all $n \in \mn$ for some constant $D_1 \in \mr _+$.
Also $\vert R_n (b_1 ) \ldots R_n (b_k ) \vert \leq   \frac{1}{n^2} D_2 $ for all $n \in \mn$ for some $D_2 \in \mr _+ $, and for a suitable constant $D_3$
\beao
 && \vert \psi (b_1 ) R_n (b_2 ) \ldots R_n (b_k )   
 + R_n (b_1 ) \psi (b_2 ) R_n (b_3 ) \ldots R_n (b_k ) + \ldots   \\
 && \qquad\qquad + R_n (b_1 ) \ldots R_n (b_{k - 1} ) \psi (b_k ) \vert \leq \frac{1}{n^2} D_3 
 \eeao
 if $k \geq 2$ so that
 $$
 S(\frac{\psi}{n} + R_n ) (b_1 \ot_s \ldots \ot _s b_k ) \leq \frac{1}{n^2} D
 $$
 if $k \geq 2$ for some $D \in \mr_+$.
 For $k = 1$ we have
 $S(\frac{\psi}{n} + R_n ) (b) = \frac{\psi}{n} (b) + R_n (b)$. Moreover, $S(\frac{\psi}{n} + R_n ) (\ei ) = 1$. It follows that
 $$
 S(\frac{\psi}{n} + R_n ) = S(0) + \frac1n D(\psi ) + \tR _n
 $$
 with $\tR _n : S(\SB ) \to \mc$ linear, $\tR _n (\ei ) = 0$, and $\vert \tR _n (x) \vert \leq \frac{1}{n^2} \tC _x$ for all $n \in \mn$ for a suitable constant $\tC _x$.
 By Lemma \ref{lemma1}
 $$
 S(\frac{\psi}{n} + R_n ) = \bigl( S(0) + \frac1n D(\psi ) + \tR _n \bigr) ^{\star n} \longrightarrow \exstar D(\psi )
 $$
 point-wise.$\square$ 
\par
We say that {\it Schoenberg correspondence holds} on a dual semigroup $\SB$ if the convolution exponential
$\exstar \psi$ is a state for each
conditionally positive linear functional $\psi$  on $\SB $ and for each positive $\odot$-independence.
Then Schoenberg correspondence holds on $\SB$ iff the CCSG of states (with respect to a positive $\odot$-independence) on $\SB$ are exactly given by $\exstar (t \psi )$ with $\psi$ conditionally positive. 
We will show that Schoenberg correspondence holds on \emph{all} dual semigroups (Theorem \ref{theorem}).
\par
Let $(\SB , \GD )$ and $(\SC , \GL )$  be two dual semigroups and let $\kappa : \SC \to \SB$ be a $*$-algebra homomorphism.
For a linear functional $\psi$ on $\SB$ we put
\begin{align*}
\gg _t &= \exstar (t \psi )  \\
\phi _t &= \exstar (t (\psi \circ \kappa ))
\end{align*}
for $t \in \mr _+$.
\bprop
\label{approx}
$$
( \gg _{t/n} \circ \kappa ) ^{\star n} \to \phi _t
$$
point-wise for all $t \in \mr _+$.
\eprop
\pn
{\it Proof}:
We can assume that $ t = 1$.
We will show that there are linear functionals $R_n : \SC \to \SB$ and constants $C_b \in \mr _+$, $b \in \SC$, such that
$$
(\gg _{1/n} \circ \kappa )(b) = \frac1n (\psi \circ \kappa )(b) + R_n (b)
$$
with
$$
\vert R_n (b) \vert \leq \frac{1}{n^2 } C_b
$$
for all $n \in \mn$.
An application of Lemma \ref{lemma2} will then prove the proposition.
\par
We have
\beao
(\gg _{1/n} \circ \kappa )(b) &=& \gg_{1/n} (\kappa (b))    \\
&=& \bigl( \exstar \frac1n D(\psi ) \bigr) (\kappa (b))  \\
&=& \frac1n D(\psi ) (\kappa (b)) + \frac{1}{2 n^2 } D(\psi ) ^{\star 2} (\kappa (b))   \\
&& \qquad + \frac{1}{3 ! n ^3 } D(\psi ) ^{\star 3}(\kappa (b)) + \ldots  \\
&=& \frac1n \psi (\kappa (b)) + \frac{1}{n^2} R_n (b)
\eeao
with $\vert R_n (b) \vert \leq C_b$ for all $n \in \mn$ for a suitable constant $C_b$.$\square$
\bprop
\label{schoenberg}
Suppose that Schoenberg correspondence holds on $\SB$.
Then for a conditionally positive linear functional $\psi$ on $\SB$,
we have that, given a positive $\odot$-independence, $\phi _t = \exstar (t (\psi \circ \kappa ))$ is a CCSG of states on $\SC$.
\eprop
\pn
{\it Proof}:
We fix a positive $\odot$-independence.
If Schoenberg correspondence holds on $\SB$, then $\gg_{t/n}$ are states on $\SB$. This implies that $\gg _{t/n} \circ \kappa$ are states on $\SC$ which by the positivity of the independence gives that $(\gg _{t/n} \circ \kappa )^{\star n}$ are states on $\SC$.
By Proposition \ref{approx} $\phi _t$ is the point-wise limit of $(\gg _{t/n} \circ \kappa )^{\star n}$.
Since the point-wise limit of states is a state, the proposition follows.$\square$
\par
Let $\psi$ be a conditionally positive linear functional on a $*$-algebra $\SA$. We form the left ideal $\SN _{\psi} = \{ a \in \SA \mit \psi (ba ) = 0 \ \forall b \in \SA \}$.
The quotient space $D = \SA / \SN _{\psi }$ is an inner product space with inner product
$\langle \eta (a) , \eta (b) \rangle = \psi (a ^*  b )$ where $\eta : \SA \to D$ denotes the canonical map.
Moreover, $\rho (a) \eta (b) = \eta (ab)$ defines a $*$-representation of $\SA$ on $D$, i.e. a $*$-algebra homomorphism from $\SA$ to the $*$-algebra $\RL (D)$ of adjointable linear operators on $D$; see \cite{Fra06,MSchue93}.
We have
\bprop
\label{determined}
If $(a_i )_{i \in I}$ is a set of generators of the algebra $\SA$, the maps $\rho$, $\eta$ and $\psi$ are determined by their values on the $a_i$.$\square$
\eprop
Since a $*$-vector  space $\SV$ generates the tensor $*$-algebra $\RT (\SV )$, we have that conditionally positive linear functionals are given by an inner product space $D$, a $*$-map $\rho : \SV \to \RL (D)$, a linear map $\eta : \SV \to D$, and a hermitian linear functional $\psi : \SV \to \mc$; cf. \cite{additive,Fra06}.
\par
Clearly, a linear functional on a $*$-algebra $\SA$ is a state if it is the expectation of a quantum random variable.
Consider the case of tensor independence. 
It is well known  that the Gelfand-Naimark-Segal (GNS) representation $\pi$ of $\phi = \exstar \psi$, $\psi$ a conditionally positive linear functional on $\RT (\SV )$, is given by
$$
\pi (v) = A^* (\eta (v)) + \GL (\rho (v)) + A (\eta (v ^* )) + \psi (v)  \id
$$
$v \in \SV$, where $A ^ * , \GL , A$ are the creation, preservation and annihilation operators on Bose Fock space over the completion $H$ of the inner product space $D$, i.e.
$$
\phi (v_1 \ot \ldots \ot v_n )  = \RE (  \pi (v_1 ) \ldots \pi (v_n ) )
$$
where expectation is taken in the vacuum state of the Fock space, see e.g. \cite{additive}.
It follows that $\phi$ is a state.
The analogous result holds in the free case if Bose Fock space is replaced by the full Fock space and $A^* , \GL, A$ are the free creation, preservation and annihilation operators; see for example \cite{GlScSp}.
In the remaining three cases of independences we can apply the reduction theory of Franz \cite{Fra06} to realize the GNS of $\phi$ on a Bose Fock space. 
We have
\bprop
\label{additive}
Schoenberg correspondence holds on tensor $*$-algebras.$\square$
\eprop
We will apply Proposition \ref{schoenberg} to the following situation.
The tensor $*$-algebra $\RT (\SB )$ of a dual semigroup $\SB$ (viewed as a $*$-vector space) carries another dual group structure than that given by the primitive comultiplication $b \mapsto i _1 (b) + i _2 (b)$.
The second dual semigroup structure is given by extending the map 
$$
(i _{\SB } \sqcup i _{\SB } ) \circ \GD : \SB \to
\RT (\SB ) \sqcup \RT (\SB ) ,
$$
with $i _{\SB} : \SB \to \RT (\SB )$ the natural embedding, to a homomorphism 
$$
\RT (\GD ): \RT (\SB ) \to \RT (\SB ) \sqcup \RT (\SB ) .
$$
We denote by $M : \RT (\SB ) \to \SB$ the multiplication map $\RT (\id )$.
\bprop
\label{alg}
We have for linear functionals $\phi _1$, $\phi _2$, $\psi$ on a dual semigroup $\SB$
\pn
{\rm (a)}
$$
(\phi _ 1 \circ M ) \star _{\RT (\GD )} (\phi _2 \circ M) = (\phi _1 \star _{\GD} \phi _2 ) \circ M
$$
{\rm (b)}
$$
\exp _{\star \RT (\GD )} (\psi \circ M) = (\exp _{\star \GD} \psi ) \circ M
$$
\pn
{\rm (c) } The linear functional $\exp _{\star \RT (\GD )} (\psi \circ M)$ vanishes on the two-sided $*$-ideal $\kernel M$.
\eprop
\pn
{\it Proof}: From (A1)
$$
(\phi _1 \circ M ) \bullet (\phi _2 \circ M) = (\phi _1 \bullet \phi _2 ) \circ (M \Sqcup M ) .
$$
Since $(M \Sqcup M ) \circ \RT (\GD ) = \GD \circ M$,
part (a) follows.
(b) is a consequence of (a) and \ref{schoenberg} for $R_n = 0$.
(c) follows from (b).$\square$
\bthm
\label{theorem}
Schoenberg correspondence holds for all dual semigroups.
\ethm
\pn
{\it Proof}:
For a positive $\odot$-independence and a conditionally positive linear functional $\psi$ on a dual semigroup $(\SE , \GD )$ we have that $\psi \circ M$ is conditionally positive on the tensor $*$-algebra $\RT (\SE )$.
Now apply Proposition \ref{schoenberg} to $\SC = (\RT (\SE ), \RT (\GD ))$,
$\SB = \RT (\SE )$ with the primitive comultiplication, and $\kappa = \id$. By Proposition \ref{additive} Schoenberg correspondence holds on $\SB$.
By Propositon \ref{schoenberg} this means that Schoenberg correspondence holds on $\SC$.
Thus $\exstar (\psi \circ M)$ is a state on $\SC$.
But, using (c) of Proposition \ref{alg} and the fact that 
$\RT (\SE ) / \ker M = \SE$, it follows that $\exstar \psi$ is a state on $\SE$.$\square$

\section{Quantum L\'evy processes}\label{sec:levy}

A {\it quantum L\'evy process} (QLP)
on a dual semigroup with respect to a positive $\odot$-independence (over a quantum probability space $(\SA , \RE )$) is a family of quantum random variables $j_{st} : \SB \to \SA$, $0 \leq s \leq t$, such that
\begin{align}
\label{QLP1}
&(j_{rs} \sqcup j_{st}) \circ \GD  = j_{rt} \ \mbox{for all} \ 0 \leq r \leq s \leq t  \\
\label{QLP2}
&j_{t_1, t_2}, \ldots , j_{t_n , t_{n + 1}} \ \mbox{are independent for all} \ n \in \mn , 0 \leq t_1 \leq \ldots \leq t_{n + 1}  \\
\label{QLP3}
&\RE \circ j_{st} \ \mbox{only depends on} \ t - s  \\
\label{QLP4}
&\lim_{t \to 0+} (\RE \circ j_{0 t} )(b) = 0 \ \mbox{for all} \ b \in \SB
\end{align}
Property \eqref{QLP1} is the increment property, \eqref{QLP2} expresses the independence of increments with respect to the underlying $\odot$-independence, \eqref{QLP3} reflects the stationarity of increments, and \eqref{QLP4} is a condition of weak continuity. 
\par
Let $\SI \subset \mr _+$ be a compact interval or equal to $\mr _+$.
Denote by $M$ the set
$$
M = \{  \gs \subset \SI  \mit 1 < \# M < \infty \}
$$
of finite subsets of $\SI$ with the inclusion of sets as partial ordering.
We write $\gs = \{ t_1 < \ldots < t_{n + 1} \}$ for a set 
$\gs = \{ t_1 , \ldots , t_{n + 1} \} \in M$, $t_1 < \ldots < t_{n + 1}$.
Define the $n$-th comultiplication $\GD : \SB \to \SB ^{\sqcup n} $, $n \in \mn _0$, of a dual semigroup $(\SB , \GD )$ recursively by
\begin{eqnarray*}
\GD _0 = 0 ;  \ \GD _{n + 1} = (\GD \Sqcup \id) \circ \GD _n
\end{eqnarray*}
For $\{ s < t \} \in M$ we denote by $\SB _{st}$ a copy of $\SB$ and by
$\iota _{st} : \SB \to \SB _{st}$ the identification map.
Put $\SB _{\gs} = \bigsqcup_{l = 1}^{n} \SB _{t_l , t_{l + 1}}$ and let
$$
f_{\{ t_1 , t_{n + 1} \} , \gs } : \SB _{t_1 , t_{n + 1}} \to \SB _{\gs}
$$
be the mapping 
$$
( \iota _{t_1 , t_2} \Sqcup \ldots \Sqcup \iota _{t_n , t_{n + 1}} ) \circ \GD _n \circ \iota _{t_1 , t_{n + 1}} ^{-1} .
$$
Moreover, for $\gs = \{ t_1 < \ldots < t_{n + 1} \}$ and $\tau \supset \gs$,
\begin{align*}
& \tau = \{ t_1 = t_{11} < \ldots < t_{1,m_1} < t_{1,m_1 + 1} = t_2 = t_{21} < \ldots < \ldots t_{2, m_2} < t_{2, m_2 + 1} \\
& \qquad = t_3 < \ldots 
< t_{n - 1 ,m_{n - 1} + 1} = t_n = t_{n1} < \ldots < t_{n, m_n} < t_{n, m_n + 1} = t_{n + 1} \} ,
\end{align*}
we define $f_{\gs \tau} : \SB _{\gs} \to \SB _{\tau}$ by
\begin{equation*}
f_{\gs \tau} = f_{\{ t_1 ,  t_2 \} ,\{ t_{11} , \ldots , t_{1,m_1 + 1} \} } \Sqcup \ldots \Sqcup
f_{\{ t_n ,  t_{n + 1} \} ,\{ t_{n1} , \ldots , t_{n,m_n + 1} \} } .
\end{equation*}
For $\tau = \{ t_1 < \ldots < t_{n + 1} \}$ and $\gs = \{ t_k < \ldots < t_l \}$, $k , l \in \{ 1, \ldots , n + 1 \}$
we put $f_{\gs \tau} : \SB _{\gs} \to \SB _{\tau}$ equal to the natural embedding.
For the general case $\tau \supset \gs$ when $\gs = \{ s_1 < \ldots < s_{m + 1} \}$ with $s_1 = t_k$, $s_{m + 1} = t_l$ for
$k , l \in \{ 1, \ldots , n + 1 \}$ we put
$$
f_{\gs \tau} = f_{\{ t_k < \ldots < t_l \} , \tau } \circ f_{ \gs , \{ t_k < \ldots < t_l \}} .
$$
Let $\psi$ be conditionally positive on $\SB$.
By Schoenberg correspondence $\phi _t = \re _{\star} ^{t \psi}$ form a convolution semigroup
of states on $\SB$.
For $\gs = \{ t_1 < \ldots < t_{n + 1}\} \in M$ we put
$$
\phi _{\gs} = (\phi _{t_2 - t_1} \bullet \ldots \bullet \phi _{t_{n + 1} - t_n} )
\circ (\iota _{t_1 , t_2} ^{-1} \Sqcup \ldots \Sqcup \iota _{t_n , t_{n + 1}} ^{-1} )
$$
to obtain a state on $\SB _{\gs}$.
The family $(\SB _{\gs}, \phi _{gs}, f _{\gs \tau} )$ is an inductive system in the category formed by
pairs $(\SB , \phi )$ where $\SB$ is a $*$-algebra and $\phi$ a state on $\SB$.
It is not difficult to see that inductive limits exist in this category.
Let $(\SA , \RE , f_{\gs })$ be the inductive limit of the above inductive system.
Then $\SA$ is a $*$-algebra, $\RE$ a state on $\SA$ and $f_{\gs } : \SB _{\gs} \to \SA$ are
$*$-algebra homomorphisms such that $\RE \circ f_{\gs} = \phi _{\gs}$ and
$f _{\tau} \circ f_{\gs \tau} = f _{\gs}$.
\bthm
\label{kolmogorov}
$$
j_{st} = f_{\{ s , t \} } \circ \iota _{st} : \SB \to \SA , \  s < t
$$
defines a QLP on $\SB$ whose convolution semigroup of states is given by $\phi _t = \re ^{t \psi } _ {\star}$.
\ethm
\pn
{\it Proof}:
We have
\begin{eqnarray*}
j_{rs} \star j_{st} &=& (j_{rs} \sqcup j_{st} ) \circ \GD   \\
&=& f_{ \{ r, s , t\} } \circ (f_{\{ r , s \} , \{ r, s , t \} } \amalg f_{\{ r, t \} , \{ r, s, t \} } ) \circ \GD   \\
&=& f_{\{ r, s , t \} } \circ f_{\{ r, t \} , \{ r, s , t \} } \\
&=& j_{rt} .
\end{eqnarray*} 
Next
\begin{eqnarray*}
\RE ( j_{t_1 , t_2 } \sqcup \ldots \sqcup j_{t_n , t_{n + 1}}) &=& \RE \circ j_{\{ t_1 , \ldots , t_{n + 1} \} }  \\
&=& \phi _{\{ t_1 , \ldots , t_{n + 1} \} }  \\
&=& \phi _{t_2 - t_1} \bullet \ldots \bullet \phi _{t_{ n + 1} - t_n}  \\
&=& (\RE \circ j_{t_1 , t_2} ) \bullet \ldots \bullet (\RE \circ j_{t_n , t_{n + 1}} )
\end{eqnarray*}
which gives the independence of increments.
Since $\RE \circ j_{st} = \phi _{t - s}$, stationarity of increments holds, too.
Weak continuity follows from the continuity of $\phi _t$. 
We have that $j_{st}$ is a QLP.
Its convolution semigroup clearly is given by $\phi _t$.$\square$
\bn
{\bf Example 1.}
A classical time-indexed stochastic process on the group $\SU _d$ of unitary $d \times d$-matrices is a process $U_t$, $t \geq 0$, of unitary operators on the Hilbert space $\mc ^d \ot \RL ^2 (\GO )$.
The process $U_t$ is a L\'evy process if the $*$-algebras $\SA _{st} \subset \RL ^{\infty} (\GO ) \subset \SB (\RL ^2 (\GO ))$, $0 \leq s \leq t$, generated by the entries $(U_{st} )_{kl}$, $k, l = 1, \ldots , d$, are independent for disjoint intervals where we write $U_{st}$ for the increment $U_s ^{-1} U_t \in \SB (\mc ^ 2 \ot \RL ^2 (\GO )) \ \cong
\RM _d (\SB (\RL ^2 (\GO ))$. This means
\begin{eqnarray*}
\RE (a_1 \cdots a_n ) = \RE (a_1 ) \cdots \RE (a_n )
\end{eqnarray*}
for $a_i \in \SA _{t_i , t_{i + 1}}$, $i = 1, \ldots , n , t_1 < \ldots < t_{n + 1}$. Here $\RE (a) = \int _{\GO} a \, \rd \RP = \langle \xi , a \, \xi \rangle$, $\xi$ the constant finction equal to 1.
Moreover, the expectation restricted to $\SA _{st}$ is only to depends on $t - s$, and $U_t$ coverges to the identity in the sense that
\begin{eqnarray}
\label{unitary}
\langle \xi , (U_{st} ^{(*)}) _{k_1 , l_1} \cdots (U_{st} ^{(*)}) _{k_n , l_n} \to \gd _{k_1, l_1 } \cdots \gd _{k_n , l_n}
\end{eqnarray}
for $t \downarrow s$. Since $\SA _{st} \subset \RL ^{\infty} (\GO )$ for all $0 \leq s \leq t$, the algebras $\SA _{t_i , t_{i + 1}}$ commute.
We pass to the noncommutative case (cf. \cite{vW84}) by considering processes $U_t$ of unitary,  on a Hilbert space $\mc \ot \SH$, $H$ a Hilbert space, with expectation given by a unit vector $\xi$ in $\SH$.
The algebras $\SA _{st}$ are defined as before, independence is in the state given by $\xi$ in the sense of a fixed noncommutative independence. Stationarity of increments is still defind by (\ref{unitary}).
\par
Let $\sku$ be the Kre\^ in dual of the compact group $\SU _d$ that is the Hopf $*$-algebra generated by commuting indeterminates $x_{kl}$ and $x_{kl} ^*$ with relations $x ^* x = \ei$ in matrix form where $x = (x_{kl})_{kl}$ and $x ^* = (x_{kl} ^* )_{kl}$.
The comultiplication is given by $\GD x = x \ot x$, the counit by $\gd x = 1$.
L\' evy processes on $\SU _d$ and QLPs on $\sku$ are in 1-1-correspondence via $j_{st} (x_{kl}) = (U_{st} ) _{kl}$.
Denote by $\SKU$ the dual group generated by non-commuting indeterminates, again given by the entries of $x$ and $x ^*$ with relations 
$x ^* x = \ei = x x ^*$ and with comultipication defined in the same manner as in the commutative case.
Notice that now $\GD$ is a map from $\SKU$ to $\SKU \sqcup _{\ei} \SKU$
that is we consider $\SKU$ as a unital dual group. For example, $\GD x_{kl} = \sum_{n = 1} ^d \iota _1 (x_{kn}) \iota _2 (x_{nl})$.
The antipode is the $*$-algebra automorphism given by $S x_{kl} = x_{lk} ^*$.
Noncommutative unitary L\'evy processes as described above and QLPs on $\SKU$ are the same objects, again via  $j_{st} (x_{kl}) = (U_{st} ) _{kl}$.
\par
Let $\psi$ be a conditionally positive hermitian linear functional on $\SKU$ and let $(D, \rho , \eta )$ be as  in Section \ref{sec:schoen}. Denote by $H$ the completion of $D$. Since the $x_{kl}$ generate $\SKU$ as a $*$-algebra, $\rho$ is determined by the unitary operator $W = \rho (x_{kl} )_{kl}$ on $\mc ^d \ot H$. 
Moreover,
\begin{eqnarray*}
0 = \eta ( \sum_{n = 1} ^d x_{nk}^* x_{nl} ) =
\sum_{n =1}^d (W^* _{nk} \eta (x_{nl}) + \eta (x_{nk} ) \gd _{nl})
\end{eqnarray*}
and $\eta (x^* _{kl} ) = - \sum_{n =1} ^d W^* _{nl} \eta (x_{nk} )$.
It follows from Proposition \ref{determined} that the QLP with generator $\psi$ is determined by $W$, the matrix $L \in \RM _d ( H)$, $L_{kl} = \eta (x_{kl})$, and the matrix $G \in \RM _d (\mc )$, $G_{kl} = \psi (x_{kl} )$.
Conversely, each such triplet $(W, L , G )$ defines a conditionally positive hermitian linear functional $\psi$ on $\SKU$ by the inductive limit construction of this section.
We described the QLPs on $\SKU$ by triplets $(W, L, G)$ consisting of a unitary on $\mc ^d \ot H$, a $d \times d$-matrix $L$ with entries in $H$, and a scalar $d \times d$-matrix $G$.
\bn
{\bf Example 2.}
Let $\mf _n$ denote the free group with $n \in \mn$ generators $g_1 , \ldots , g_n$.
The group algebra $\mc\mf_n$ of $\fn$ is a dual group with $g ^* = g ^{-1}$, $\GD : \cfn \to \cfn \sqcup _{\ei} \cfn$, $\GD g_i =
\iota _1 (g_i ) \iota _2 (g_i )$, $\gd g_i = 1$.
A QLP on $\fn$ (that is on $\cfn$) is given by a vector $(U_t ^{(1)}, \ldots  \ldots , U_t ^{(n)} )$ of unitary operators $U_t ^{(i)} $, $i = 1, \ldots , n$, $t \in \mr _+$, on a Hilbert space $\SH$ such that the $*$-algebras $\SA _{st} \subset \SB (\SH )$ generated by the operators $(U_s ^{(i)} )^{-1} U_t ^{(i)}$ are independent for disjoint intervals in a state on $\SB (\SH )$ given by a unit vector in $\SH$.
Stationarity of increments and continuity are defined as in the case of $\SKU$. We find that conditionally positive hermitian linear functionals, and thus QLPs, on $\cfn$ are given by triplets $(W , L , G)$ consisting of a vector $W = (W ^{(1)}, \ldots , W^{(n)} )$ of unitary operators on a Hilbert space $H$, a vector $L = (L^{(1)}, \ldots , L^{(n)} )$ of elements of $H$ and $G \in \mc ^n$. 
\bn
\pn

\bibliographystyle{plain}
\bibliography{mybib-1}

\enddocument